\magnification=1200
\input amssym.def
\input epsf

% Math:
\def \BZ {{\Bbb Z}}

% Text:

\def\sqr#1#2{{\vcenter{\vbox{\hrule height.#2pt
     \hbox{\vrule width.#2pt height#1pt \kern#1pt
     \vrule width.#2pt} \hrule height.#2pt}}}}
\def\square{\ $\sqr66$}

\centerline {\bf A MOVE ON DIAGRAMS THAT GENERATES S-EQUIVALENCE OF KNOTS}
\bigskip
\centerline {Swatee Naik and Theodore Stanford \footnote \dag 
{Research supported in part by the Naval Academy Research Council}}
\bigskip

\bigskip
\noindent
Abstract: \ Two knots in three-space are S-equivalent 
if they are indistinguishable
by Seifert matrices.  We show that S-equivalence is
generated by the doubled-delta move on knot diagrams.
It follows as a corollary that a knot has trivial Alexander
polynomial if and only if it can be undone by 
doubled-delta moves.  

\bigskip
\bigskip

We consider tame, oriented knots in oriented $S^3$, with
equivalence being ambient isotopy.  A Seifert surface
for such a knot is an oriented surface whose boundary
is the given knot, and whose orientation induces the given
orientation on the knot.  An oriented surface $S$ in 
$S^3$ has a linking form $\langle *,* \rangle$
on the homology $H_1 (S)$, where $\langle x,y \rangle$ is defined
to be the linking number of the cycle $x$ with
the cycle $y$ slightly pushed off $S$ in a direction
determined by the orientation of $S$.
Given a knot $K$,
choose a Seifert surface $S$ for $K$ and a basis
for $H_1 (S)$.  Then the linking form is represented
by an integer matrix $M$, which is called a Seifert
matrix for $K$.  Two knots are called $S$-equivalent
if they have a common Seifert matrix (which is the same
as saying that they have a common Seifert form).
We sketch a proof at the end of the paper that this
is an equivalence relation.  The reader who wishes
may ignore this proof and take S-equivalence to
be the smallest equivalence relation that includes
any pair of knots which have a common Seifert matrix.
Two knots $K$ and $K^\prime$ are then S-equivalent
if and only if
there exists a sequence $K = K_1, K_2, \dots K_m = K^\prime$,
such that for all $1 \le i < m$,
$K_i$ and $K_{i+1}$ have a common Seifert matrix.
It makes no difference in the proof of 
Theorem~A which definition we use, 
since in either case what we need to show
is that two knots share a common Seifert matrix if and only
if they are equivalent by 
certain diagram moves which we will call
doubled-delta moves.

The usual way to define S-equivalence is to define it
first for matrices, and then to define it for knots by
saying that knots with S-equivalent matrices are S-equivalent.
See Gordon~[3]
and Kawauchi~[5]
for the standard definition, for further references,
and for more detail on the following statements.
S-equivalence of matrices was first introduced by Trotter in 
[13] under the name 
$h$-equivalence.
Murasugi [10]
and Rice [11]
applied it to 
matrices obtained from knot diagrams. 
None of the abelian invariants, such
as the Alexander polynomials, homology of cyclic and branched covers,
or signatures, can distinguish between S-equivalent knots.
It was shown by 
Levine [7]
that in higher dimensions simple knots
are characterized by S-equivalence. Two knots are 
S-equivalent if and only if their (integral) Blanchfield pairings are 
isometric.  This follows from work
of Levine~[7]
and Kearton~[6],
and was also proved by 
Trotter~[14]
from a purely algebraic point of view.

A knot may be given by a regular projection in the usual
way, with equivalence of diagrams given by the Reidemeister
moves.  For more details on knots, diagrams, and
Seifert surfaces and matrices, see 
Rolfsen~[12].
or Kawauchi~[5].
We consider now the delta move and the
{\it doubled-delta move},
shown in Figure~1.

\vfil
\eject

\epsfysize = 5truecm
\centerline {\epsffile {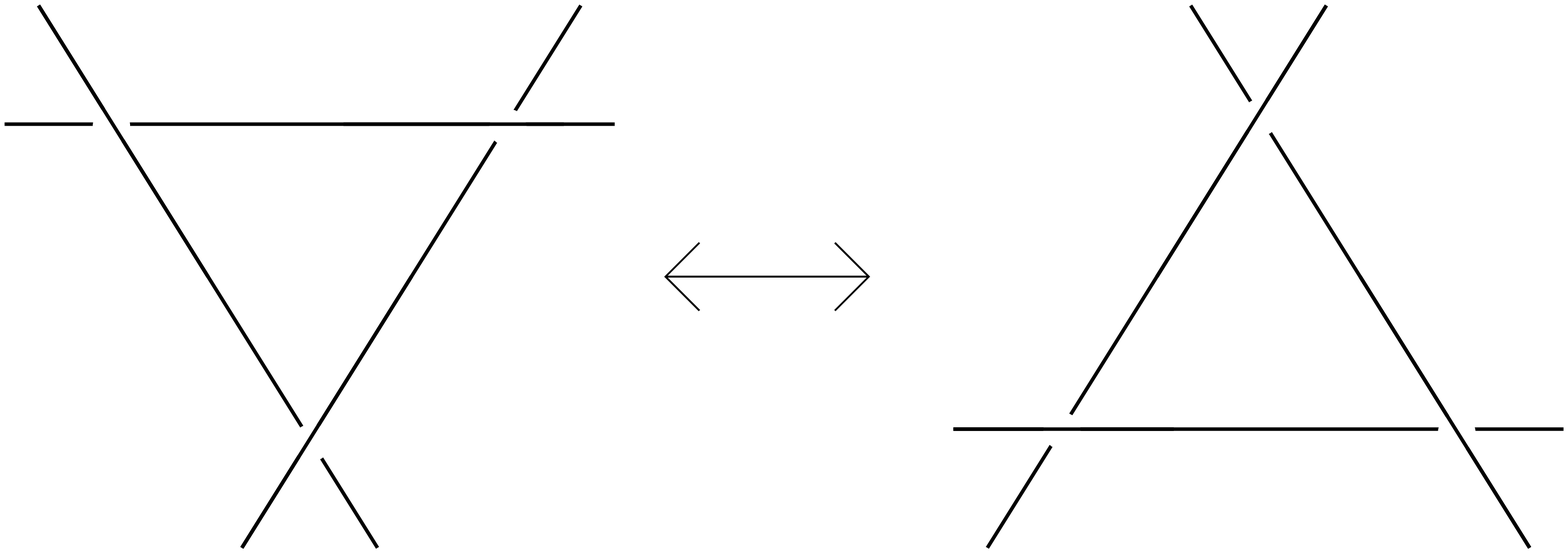}}
\medskip
\centerline 
{Figure 1a. The delta move.}

\bigskip
\epsfysize = 5truecm
\centerline {\epsffile {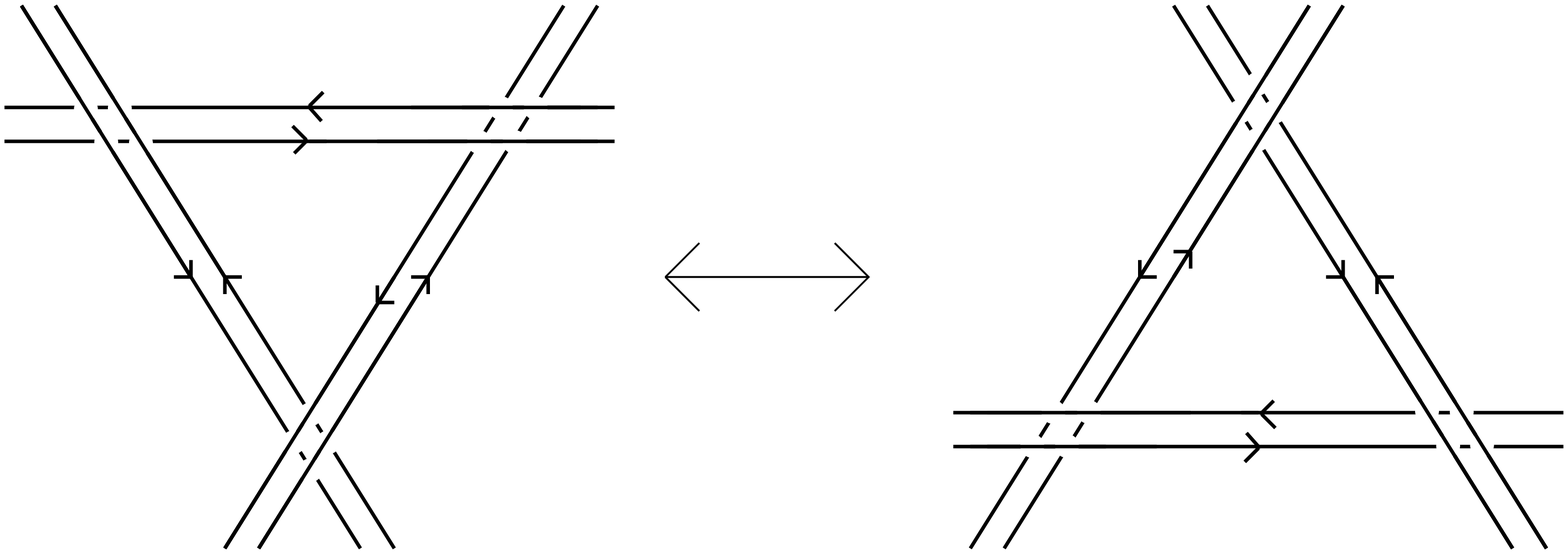}}
\medskip
\centerline 
{Figure 1b. The doubled-delta move.}
\bigskip

The orientations shown on the strands 
of the doubled-delta move
matter only in that
the arrows on each pair of parallel 
strands are oriented oppositely.
Any two such choices of orientation give essentially the
same move, since cancelling 
half-twists may be added to pairs
of strands just outside the disk where the move takes place.

\bigskip
\noindent
{\bf Theorem A.} \ 
Two knots $K$ and $K^\prime$
are S-equivalent if and only if they
are equivalent by a sequence of 
doubled-delta moves.
\medskip

Any knot with a trivial Alexander polynomial has a
trivial Alexander module and therefore a trivial
Blanchfield pairing.  Therefore, any two knots
with trivial Alexander polynomial are S-equivalent.
(Note, however, that it is not true in general
that knots with the same Alexander polynomial
are S-equivalent.)
Thus we have

\medskip
\noindent
{\bf Corollary B.} \  
A knot may be undone by doubled-delta
moves if and only if it has trivial Alexander polynomial.
\medskip

Since the doubled-delta move takes place inside a disk,
S-equivalence and S-triviality may be considered to
be local properties in some sense.

Theorem~A and Corollary~B are reminiscent of the result
of Kauffman~[4], 
which states that two knots have the
same Arf invariant if and only if they are equivalent
by a sequence of band-pass moves as in Figure~2.  The
idea is that a band-pass move can be used to undo
the knotting and linking of the bands in a Seifert surface,
with the Arf invariant (which is always either 0 or 1) being
the only obstruction to completely trivializing the
surface and thus the knot.  The doubled-delta move
is used to undo the knotting and linking of the bands
of a Seifert surface in a similar way.  However, since
the delta move preserves linking numbers, the doubled-delta
move cannot change the Seifert form of a surface or
of a knot.

\bigskip
\epsfysize = 2.5truecm
\centerline {\epsffile {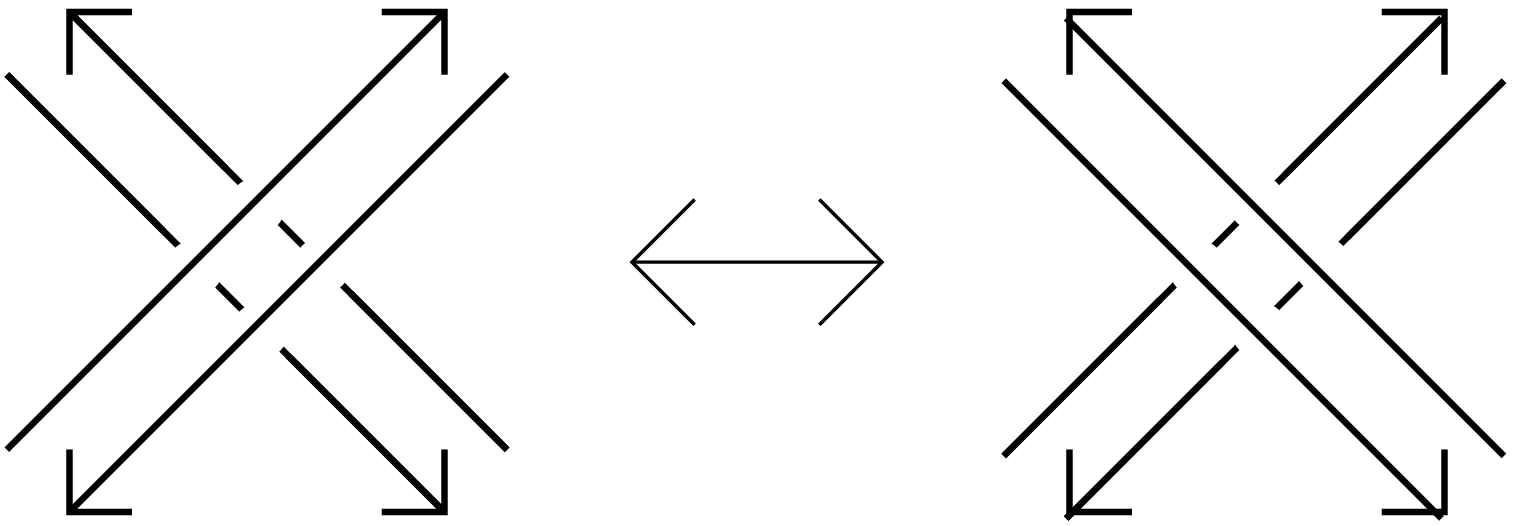}}
\medskip
\centerline 
{Figure 2. The band-pass move.}
\bigskip

Matveev~[8]
and Murakami and 
Nakanishi~[9]
have shown that two links are equivalent by delta moves if and
only if they have the same pairwise linking numbers.  We shall
need to generalize this result to string links.  Although the
proofs of [8] and
[9] appear to generalize,
we shall give a different
proof.

\medskip
\noindent
{\bf Theorem C:} \ 
Two string links are equivalent by a sequence of delta moves
if and only if they have the same pairwise linking numbers.
\medskip

\noindent
{\it Proof of Theorem A:} \ 
First we need to show that a
doubled-delta move does not change the
S-equivalence class of a knot.  Consider
a knot diagram with the left-hand side
of a doubled-delta move inside a planar disk D.
Temporarily cut the bands of the 
doubled-delta move---replace the left-hand
side of Figure~1a with Figure~3.  Apply
Seifert's algorithm to the resulting link
to obtain an oriented surface, and then
add the bands back in to obtain a Seifert
surface for the original knot. For this
surface, it is clear that the doubled-delta
move may be applied without changing the
linking form.

\bigskip
\epsfysize = 5truecm
\centerline {\epsffile {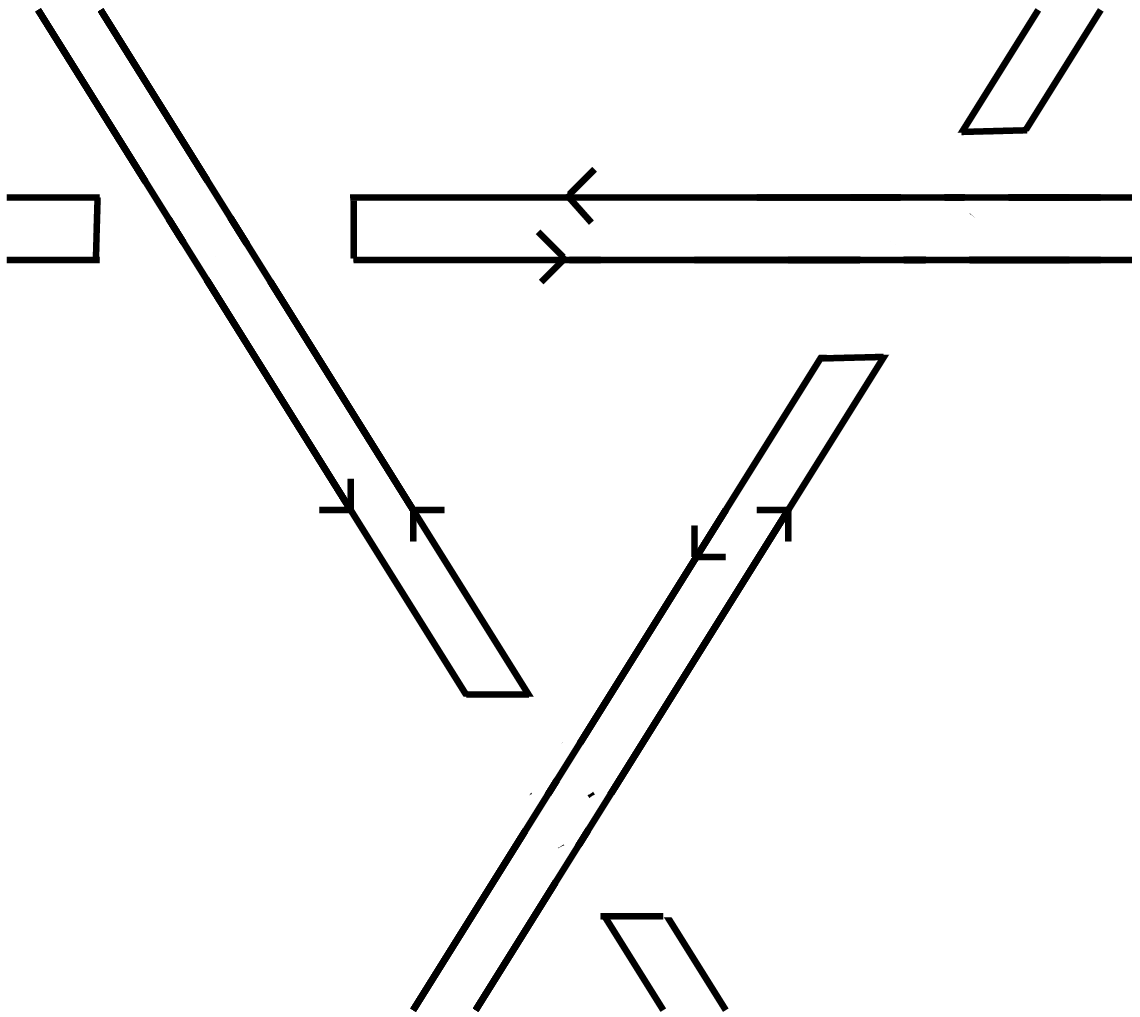}}
\medskip
\centerline 
{Figure 3.}
\bigskip

Now let $K_1$ be a knot with Seifert surface $S_1$ of genus $n$
and associated Seifert
matrix $M$, and let $K_2$ be a knot with Seifert surface $S_2$ and
Seifert matrix $M$.  Let $F_n$ be
the abstract surface with genus $n$ and
with one boundary component,
specifically realized as a disk with
bands
as shown in Figure~4.  
An ordered  basis for $H_1 (F_n)$ is
given by $\beta = ([a_1],[a_2], \dots [a_{2n}])$.
We take $F_n$ and the $a_i$
to be oriented such that 
$\langle [a_1], [a_2] \rangle = 1$ 
and $\langle [a_2], [a_1] \rangle = 0$ 
The matrix $M$ represents the linking
form of $S_1$ with respect to some
basis $B_1$ of $H_1 (S_1)$.  Let
$\phi_1: F_n \to S_1$ be an orientation-preserving
homeomorphism. Then 
$\phi_1 (\beta) 
= (\phi_1([a_1], \phi_1([a_2]), \dots \phi_1 ([a_{2n}])$ 
is also a basis
for $H_1 (S_1)$, and so there exists
an invertible integer matrix $A_1$
such that 
$N_1 = A_1 M A_1^T$ represents the
linking form of $S_1$ with respect to
the basis $\phi_1 (\beta)$.  (If $A$ is a matrix,
we denote its transpose by $A^T$ and the inverse
of $A^T$ by $A^{-T}$.) Define
$B_2, \phi_2, A_2$, and $N_2$ the same
way for the surface $S_2$.

\bigskip
\epsfysize = 5truecm
\centerline {\epsffile {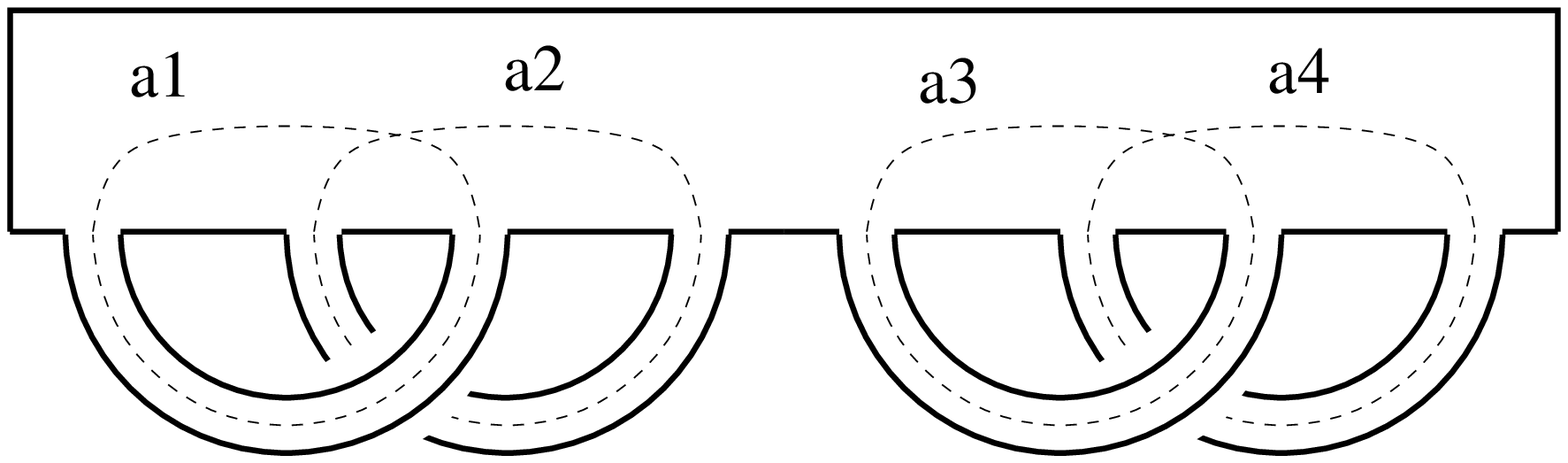}}
\medskip
\centerline 
{Figure 4. The standard surface $F_2$.}
\bigskip

For $i = 1,2$, we have $N_i - N_i^T = X_n$, where
$\displaystyle X_1 = \matrix {
0 & 1 \cr -1 & 0}$ and $X_n$ is built up with $n$ copies
of $X_1$ down the diagonal and $0$ everywhere else.
Since
$N_1 = A_1 A_2^{-1} N_2 A_2^{-T} A_1^T$,
we get $X_n = A_1 A_2^{-1} X_n A_2^{-T} A_1^{-1}$, so
$A_1 A_2^{-1} = C = [C_{i,j}]$ is
a symplectic matrix.  The mapping class group of
a surface with 0 or 1 boundary components acts on
$H_1$ of the surface, and if a basis is chosen
for $H_1$ with the same intersection properties
as the $a_i$ curves on $F_n$, then the matrices
that represent this action are symplectic matrices.
Moreover, this map from the mapping class group
to the symplectic group is well-known to be
surjective, and therefore
there exists a 
homeomorphism $g: S_2 \to S_2$ 
such that 
$[g(\phi_2(a_i))] = \sum_{j=1}^{2n} C_{i,j} [\phi_2 (a_j)]$
for all $1 \le i \le 2n$, and $N_1$ represents the
linking form of $S_2$ with respect to the ordered basis
$([g(\phi_2(a_1))], [(\phi_2(a_2))], \dots [g(\phi_2(a_{2n}))])$.
Now we may use $\phi_1$ and $g \circ \phi_2$ to put
$S_1$ and $S_2$, respectively, into a disk and band form
as in Figure~5.
The only difference now between $S_1$ and $S_2$ is
in the (framed) 
string links $L_1$ and $L_2$, and the
pairwise linking numbers of these two string
links are identical to each other, both being given
by the matrix $N_1$.  Hence by Theorem~C
there exists a sequence of delta moves taking
$L_1$ to $L_2$, and therefore a sequence
of doubled-delta moves taking $K_1$ to $K_2$.

\vfil
\eject
\epsfysize = 5truecm
\centerline {\epsffile {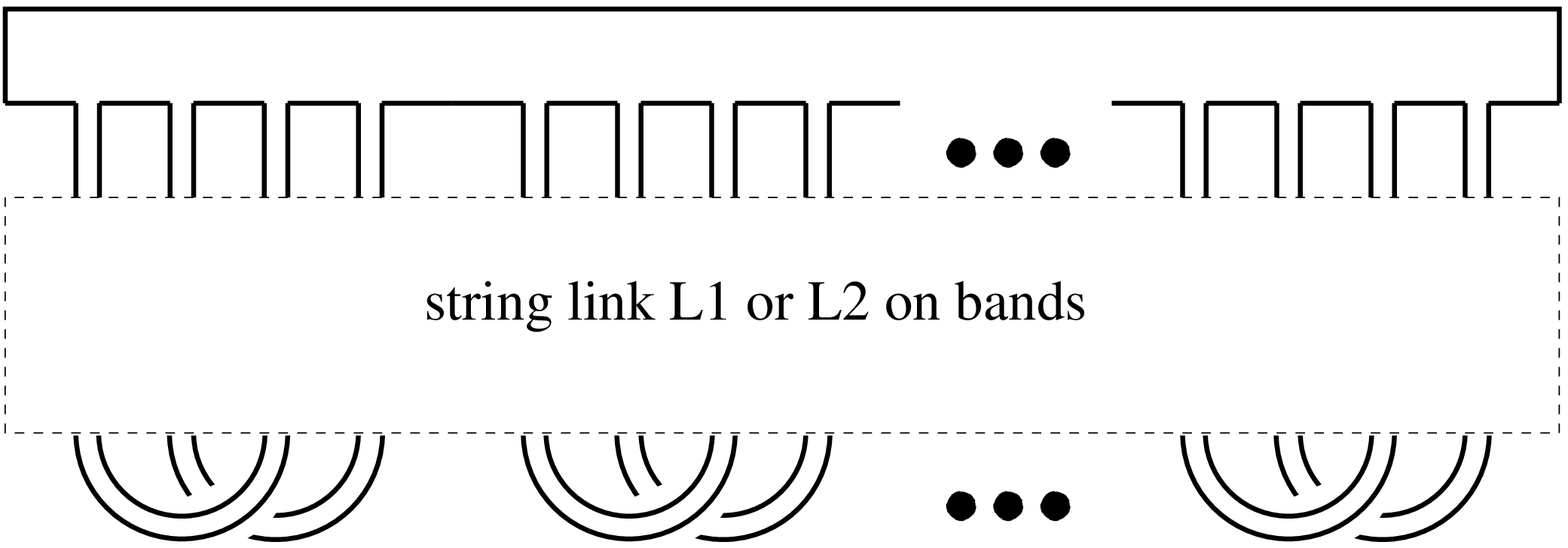}}
\medskip
\centerline 
{Figure 5. The standard form for $K_1$ and $K_2$.}
\bigskip

Note that although Theorem~C
doesn't address the framing issue, it
is clear that a delta move doesn't change
the framing on any strand of a string link,
and a doubled-delta move doesn't change the
self-linking of any of the bands in 
the surfaces $S_1$ and $S_2$.
Moreover, each band in $L_1$ must have the same
framing as the corresponding band in $L_2$, because
the framings in both diagrams are given by the diagonal
entries of $N_1$.
\square

\medskip
\noindent
{\it Proof of Theorem C:} \ 
The ``only if'' part is easy to check, so it
is left to show that two string links with matching linking
numbers are connected by a sequence of delta moves.

Let $P_n$ be the group of pure braids, those
elements of the braid group $B_n$ which induce
the identity permutation on the endpoints of
the strands.  For details and presentations of
$P_n$ and $B_n$, see Birman~[2].
For any $1 \le i<j \le n$, the
map $P_n \to \BZ$ 
which measures the
linking number between the $i$th and $j$th strands
in a pure braid is a group homomorphism.  
In fact, these are all the
abelianizing homomorphisms of $P_n$, so that
$p \in P_n$ has all its linking numbers $0$ if and
only if $p \in P_n^\prime$, the commutator subgroup
of $P_n$.
Three other facts are also easy to check.
First, any
commutator of the form 
$$p_{i,j} p_{j,k} p_{i,j}^{-1} p_{j,k}^{-1} \in P_n
\leqno {(D)}$$
may be
undone by a delta move, where $p_{i,j}$ is the standard
pure braid generator which links the $i$th and the $j$th
strands. Second, adding the commutator
(D) to a presentation of $P_n$
abelianizes the three-strand subgroup of $P_n$
generated by $p_{i,j}, p_{j,k}, p_{i,k}$.  
Third, abelianizing all such three-strand
subgroups of $P_n$ abelianizes $P_n$.  
Thus if a pure braid $p$ has all its linking numbers
$0$ then it can be undone by delta moves.

\vfil
\eject
\epsfysize = 5truecm
\centerline {\epsffile {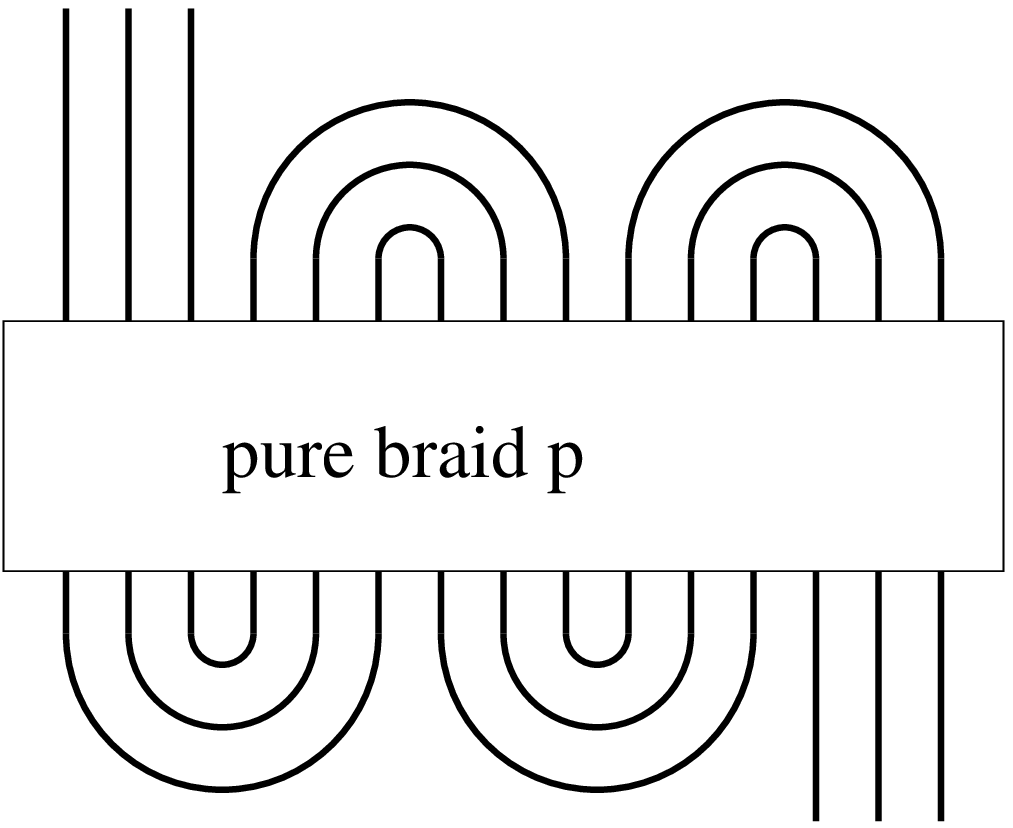}}
\medskip
\centerline 
{Figure 6. }
\bigskip

Now let $L$ be an $n$-strand string link with all pairwise
linking numbers equal to $0$.
For some $k>0$, $L$ is represented
by a diagram as in Figure~6, where $p$ is a pure braid
on $kn$ strands. Let us label the strands of $p$ with a double index
$(i,a)$, indicating the $a$th braid strand of the
$i$th string link strand.  For example, the $4$th
braid strand in Figure~6 (in the usual sense) will
be labeled $(3,2)$, and the $5$th will be labeled $(2,2)$.
Let $p_{(i,a)(j,b)}$ be the standard braid generator, as 
above, linking the $(i,a)$ strand with the $(j,b)$ strand.
We have $p_{(i,a)(j,b)}=p_{(j,b)(i,a)}$, and for notational
convenience we set $p_{(i,a)(i,a)} = 1$ for all $i,a$.
Let ${\rm lk} ((i,a),(j,b))$ be the linking number between the
$(i,a)$ strand and the $(j,b)$ strand.  We measure
this linking number with respect to the braid orientation
on the strands, which coincides with the string link
orientation exactly when the second index of the strand
is odd.
The linking number 0 condition on $L$ becomes
$$\sum_{a=1}^k \sum_{b=1}^k  
(-1)^{a+b} {\rm lk} ((i,a),(j,b)) = 0 \leqno {(E)}$$
for all $1 \le i \ne j \le n$. 
Let $b$ be even, 
$1 \le a,b \le k$, and
$1 \le i, j \le n$.
Then we may replace $p$ with
$(p_{(i,a)(j,b)}p_{(i,a)(j,b+1)})^{\pm 1}p$  
without affecting $L$.
Similarly, if $b$ is odd, then we may replace $p$ by
$p (p_{(i,a)(j,b)}p_{(i,a)(j,b+1)})^{\pm 1}$  
without affecting $L$.
The effect of multiplication by 
$(p_{(i,a)(j,b)}p_{(i,a),(j,b+1)})^{\pm 1}$
on the linking numbers of $p$ is the same in either case.
If $i=j$ and $a=b$ or $a=b+1$, then
${\rm lk} ((i,b)(i,b+1))$ goes up by one or down by one.
Otherwise, 
${\rm lk} ((i,a),(j,b))$ and ${\rm lk} ((i,a),(j,b+1))$
both either go up by one or go down by one.
Repeated multiplications by 
$(p_{(i,a)(j,b)}p_{(i,a),(j,b+1)})^{\pm 1}$
for various appropriate values of $i,j,a,b$
thus suffice to make ${\rm lk} ((i,a),(b,j)) = 0$
for all $a>1$ and $b>1$, and condition (E) then forces
all the remaining linking numbers to be 0 as well.

Finally, if $L$ and $L^\prime$ are two string links with
the same pairwise linking numbers, let $p \in P_n$ have
the same pairwise linking numbers as $L$ and $L^\prime$.
Clearly $L$ and $L^\prime$ are equivalent by delta moves
if and only if $p^{-1} L$ and $p^{-1} L^\prime$ are equivalent
by delta moves, but these last two string links have all
pairwise linking numbers $0$, and so are both equivalent
to the unlink by delta moves.
\square

\vfil
\eject
\noindent
{\it Sketch of proof that $S$-equivalence is an equivalence relation:} \ 
The reflexive and symmetric properties are obvious, so we only
need to consider transitivity.
If $M$ and $M^\prime$ are Seifert matrices, then we may say
that $M^\prime>M$ if $M^\prime$ can be obtained from $M$ by a
sequence of unimodular congruences, column enlargements, and
row enlargements.
See for example Kawauchi~[5].
These matrix operations
correspond to changing the basis of $H_1$ of the Seifert surface,
and to adding a tube to a Seifert surface.  It is well-known
that if $<>$ is the equivalence relation generated by $<$, and
$M$ and $M^\prime$ are two Seifert matrices of the same knot,
then $M<>M^\prime$.  For a recent elementary proof that 
two Seifert surfaces of a knot are tube-equivalent, see
Bar-Natan, Fulman, and 
Kauffman~[1].

We make two claims:  First, that if $M_1$ and $M_2$ are
Seifert matrices and $M_1 <> M_2$, then there exists a Seifert
matrix $M_3$ such that $M_3>M_1$ and $M_3>M_2$.  Second,
if $M_1$ is a Seifert matrix for a knot $K$, and $M_2 > M_1$,
then $M_2$ is also a Seifert matrix for $K$.  Both of these
are verified with elementary matrix and tube operations.

Now suppose that knots $K_1$ and $K_2$ share a common matrix
$M_{12}$ and that knots $K_2$ and $K_3$ share a common matrix
$M_{23}$.  Then by the first claim
there exists a Seifert matrix $M_2$ for $K_2$
such that $M_2 > M_{12}$ and $M_2 > M_{23}$, and 
by the second claim 
$M_2$ is a Seifert matrix for both $K_1$ and $K_3$.
\square

\vfil
\eject
\noindent
{\bf References.}

\smallskip
\item {[1]}
D.~Bar-Natan, J.~Fulman, and L.H.~Kauffman.
{\it An elementary proof that all spanning surfaces of a link
are tube-equivalent.}
Journal of Knot Theory and its Ramifications 7 (1998), 873--879.

\smallskip
\item {[2]}
J.S. Birman.
``Braids, Links and Mapping Class Groups.''
Annals of Mathematics Studies~82.
Princeton University Press, 1975.

\smallskip
\item {[3]}
C.McA. Gordon.
{\it Some aspects of classical knot theory},
``Knot Theory'' (ed. J. C. Hausmann), 1--60,
Lecture Notes in Mathematics 685,
Springer--Verlag, 1978.

\smallskip
\item {[4]}
L.H. Kauffman.
``Formal knot theory.''
Mathematical Notes~30.
Princeton University Press, 1983.

\smallskip
\item {[5]}
A. Kawauchi.  ``A Survey of Knot Theory.''
Birkhauser, 1996.

\smallskip
\item {[6]}
C. Kearton.
{\it Blanchfield duality and simple knots}.
Transactions of the American Mathematical Society
202 (1975), 141--160.

\smallskip
\item {[7]}
J. Levine. 
{\it An Algebraic Classification of some Knots of Codimension 2.}
Commentarii Math.~Helv., 15 (1970) 185--198.

\smallskip
\item {[8]}
S.V. Matveev.
{\it Generalized surgeries of three-dimensional manifolds and
representations of homology spheres.}
Matematicheskie Zametki 42 (1987), no. 2, 268--278.
(English translation in 
Mathematical Notes 42 (1987) 651--656.)

\smallskip
\item {[9]}
Murakami and Nakanishi.
{\it On a certain move generating link homology.}
Mathematische Annalen 284 (1989), 75--89.

\smallskip
\item {[10]}
K. Murasugi.
{\it On a Certain Numerical Invariant of Link Types},
Transactions of the American Mathematical Society 117
(1965) 387--422.

\smallskip
\item{[11]}
P.M. Rice. 
{\it Equivalence of Alexander Matrices}, 
Mathematische Annalen 193 (1971) 65--75.

\smallskip
\item {[12]}
D. Rolfsen.
``Knots and Links''
Mathematics Lecture Series 7.
Publish or Perish, Wilmington, DE, 1976.

\smallskip
\item {[13]}
H.F. Trotter. 
{\it Homology of Group Systems with applications
to Knot Theory}. 
Annals of Mathematics 76 (1962) 464--498.

\smallskip
\item {[14]}
H.F. Trotter. 
{\it On S-Equivalence of Seifert Matrices},
Inventiones Mathematicae 20 (1973) 173--207.

\bigskip
\bigskip

Department of Mathematics

University of Nevada

Reno, NV \ 89557

{\tt naik@unr.edu}

\medskip

Department of Mathematics

United States Naval Academy

Annapolis, MD \ 21402

{\tt stanford@usna.edu}

\end